\documentclass[11pt]{article}

\usepackage{amssymb,latexsym,amsmath}		% AMS-Formelkram
\usepackage{xypic}				% Diagramme
\usepackage{a4wide}				% Papierformat
\usepackage[ansinew]{inputenc}			% dt. Zeichensatz direkt
\usepackage{titlesec}

\numberwithin{equation}{section}

\titleformat{\subsection}
	{\normalfont\large\bfseries}{\thesubsection}{0.5em}{}

\usepackage[thmmarks,thref]{ntheorem}

\theoremstyle{plain}

\theorembodyfont{\normalfont\selectfont}
\theoremseparator{.}

\newtheorem{theorem}{Theorem}[section]		% Sätze etc. mit Nummer
\newtheorem{corollary}[theorem]{Corollary}
\newtheorem{lemma}[theorem]{Lemma}
\newtheorem{proposition}[theorem]{Proposition}

\theorembodyfont{\normalfont\selectfont}	% Def./Bsp./Bem. ohne
\theoremstyle{nonumberplain}
\newtheorem{definition}{Definition}
\newtheorem{remark}{Remark}
\newtheorem{example}{Example}

\theoremstyle{nonumberplain}			% Beweise
\theoremheaderfont{\itshape}
\theoremseparator{.}
\theoremsymbol{\rule{1ex}{1ex}}
\newtheorem{proof}{Proof}
\newcommand{\F}{{\mathbb{F}}}
\newcommand{\C}{{\mathbb{C}}}

\newcommand{\Z}{{\mathbb{Z}}}

\newcommand{\LL}{{\mathbb{L}}}

\newcommand{\OO}{{\mathcal{O}}}

\newcommand{\x}{x}                 % geändert!
\newcommand{\s}{{\bf s}}

\newcommand{\ph}{{\varphi}}
\newcommand{\eps}{{\varepsilon}}

\newcommand{\onto}{{\twoheadrightarrow}}

\newcommand{\Ext}{{\mathrm{Ext}}}
\newcommand{\Hom}{{\mathrm{Hom}}}
\newcommand{\Der}{{\mathrm{Der}}}

\newcommand{\mf}[1]{\mathfrak{#1}}
\newcommand{\mc}[1]{\mathcal{#1}}

\newcommand{\p}[1]{\partial #1}

\newcommand{\Singular}{\textsc{Singular}}

\newcommand{\ARA}{\begin{array}}
\newcommand{\ARE}{\end{array}}
\newcommand{\SMA}{\begin{smallmatrix}}
\newcommand{\SME}{\end{smallmatrix}}
\newcommand{\BMA}{\begin{matrix}}
\newcommand{\BME}{\end{matrix}}

\begin{document}

\pagenumbering{arabic}

%-- Titel und Abstract --------------------------------------------------------

\title{Deformations with section: Cotangent cohomology, flatness conditions and modular subgerms}

\author{
\begin{tabular}{cc}
Tobias Hirsch\thanks{Partially supported by G. I. F.} &
Bernd Martin$^*$\\
{\small email: \texttt{hirsch@math.tu-cottbus.de}} &
{\small email: \texttt{martin@math.tu-cottbus.de}}\\
{\small BTU Cottbus, Institut für Mathematik,} &
{\small BTU Cottbus, Institut für Mathematik,}\\
{\small Postfach 10 13 44, D-03046 Cottbus} &
{\small Postfach 10 13 44, D-03046 Cottbus}
\end{tabular}
}

\maketitle

\begin{abstract}
We study modular subspaces corresponding to two deformation functors associated to an isolated singularity $X_0$: the functor $Def_{X_0}$ of deformations of $X_0$ and the functor $Def^s_{X_0}$ of deformations with section of $X_0$. After recalling some standard facts on the cotangent cohomology of analytic algebras and the general theory of deformations with section, we give several criteria for modularity in terms of the relative cotangent cohomology modules of a deformation. In particular it is shown that the modular strata for the functors $Def_{X_0}$ and $Def^s_{X_0}$ of quasihomogeneous complete intersection singularities coincide. Flatness conditions for the first cotangent cohomology modules of the deformation functors under consideration are then compared.
\end{abstract}

\tableofcontents

\newpage

%-- Introduction --------------------------------------------------------------

\addcontentsline{toc}{section}{Introduction}
\section*{Introduction}

When trying to construct moduli spaces for analytic objects, one major difficulty results from the fact that, in general, semi-universal deformations of the objects under consideration are not universal. One approach to construct local moduli is the study of some kind of maximal universal locus in the base space of a semi-universal deformation. The corresponding notion of a modular deformation has been introduced by V. P. Palamodov in \cite{Pal78} for complex spaces, later on in a formal context by O. A. Laudal in \cite{Lau79}.

In these notes we study the strictly local situation of germs of analytic spaces, i. e. deformations of singularities. This has been done before, see e. g. \cite{Mar02} and \cite{Mar03}. One obtains the following characterization of modular subspaces in the base $S$ of a semi-universal deformation $\xi: X\to S$ of a singularity $X_0$: $M\subseteq S$ is modular if and only if all vector fields of the special fibre can be lifted to vertical vector fields of the family $\xi$ -- a criterion that already appeared in \cite{Pal78} in the context of compact complex spaces. For deformations of complete intersections and space curves, modularity of a subspace $M\subseteq S$ can also be interpreted as flattening stratum of the relative Tjurina-module $T^1(X/S)$ of the family. This can be used to actually compute non-trivial examples of modular strata using a new algorithm to determine local flattenings. The algorithm and its implementation in the computer algebra system \Singular, cf. \cite{GPS01}, are also explained in \cite{Mar02}.

This article investigates the similar situation of deformations with section of $X_0$, i. e. deformations $\xi:X\to S$ of $X_0$ together with a morphism $\sigma:S\to X$ such that $\xi\circ\sigma=id_S$. The first section is a brief review of the cotangent cohomology of analytic algebras and their morphisms which we will use extensively later on. 

Section 2 collects the basic facts on deformations with section. In particular we state the results in Buchweitz' thesis \cite{Buch81} on the cotangent cohomology of such deformations and on the construction of a semi-universal deformation with section, it is given explicitly for the case of complete intersections.

In section 3 we derive a variant of the Kodaira-Spencer sequence of the deformation $\xi$ that also contains information on the given section. With this tool at hand we formulate and prove a criterion of modularity as liftability of vector fields on $X_0$ with values in the maximal ideal of $\OO_{X_0}$ to vertical vector fields of the family. As an application we prove that, for quasihomogeneous complete intersections, the modular strata with respect to both deformation functors coincide, this being mainly a consequence of A. G. Alexandrov's description of the module of derivations $\Der_\C(\OO_{X_0})$ for this class of singularities in \cite{Alex85}.

In the closing section we give an interpretation of this criterion as flatness of the first cotangent cohomology module of $X$ over $S$ with coefficients in the kernel of $\sigma^*:\OO_X\to\OO_S$. We then compare flatness conditions for the zeroth and first relative cotangent cohomology modules. We finish by investigating the relationship between modular subspaces for both deformation functors and prove that both spaces coincide for deformations with singular section of hypersurfaces.

%-- Summary of Cotangent Cohomology--------------------------------------------

\section{A Summary of cotangent cohomology}

In this introductory section we summarize the facts about the cotangent cohomology of a morphism $(X,0)\to (S,0)$ of complex space germs or, equivalently, of morphisms $\mc{O}_{S,0}\to\mc{O}_{X,0}$ of analytic algebras that we will need later on. Since we are working exclusively in the category $(Gan)$ of germs of complex spaces (resp. the category $(Analg)$ of analytic algebras), we usually omit the distinguished point and simply write $X$ for a germ $(X,0)$ of a complex space.

\subsection{Cotangent cohomology of analytic algebras}

We start by recalling the construction of the cotangent complex of a morphism $A\to B$ of analytic algebras and the definition and properties of the so-called cotangent cohomology functors $T^i(B/A,-)$. This is the complex-analytic analogue of Andre-Quillen cohomology (cf. \cite{Quillen70}, \cite{Andre74}). Details and proofs of the statements below can be found in  \cite{Pal76}, \cite{Flenn78}, \cite{Buch81}, \cite{Pal82}, for example.

\begin{definition}
Let $A\to B$ be a homomorphism of analytic algebras. 
\begin{enumerate}
\item A {\em resolvent} $R$ for $B$ over $A$ is a DG-algebra (a differential graded anticommutative analytic algebra) that is a free $A$-algebra, together with a surjective homomorphism $p:R\onto B$ over $A$ which is a quasiisomorphism of complexes, i. e. $H^n(R)=0$ for $n<0$ and $H^0(R)\simeq B$, and this isomorphism is induced by $p$.
\item The complex 
\[\LL_{B/A}^\bullet:=\Omega_{R/A}\otimes_R B\]
 is called the {\em cotangent complex} of $B$ over $A$. Here $\Omega_{R/A}$ denotes the module of {\em Kähler differentials} of $R$ over $A$.
\item For a $B$-module $M$ we define the {\em cotangent cohomology modules} of $B$ over $A$ with values in $M$ as
\[T^i(B/A,M):=\Ext_B^i(\LL_{B/A}^\bullet,M),\qquad i\ge 0.\]
\end{enumerate}
\end{definition}

\begin{remark}
For the existence of a resolvent of $B$ over $A$ and the well-definedness of the functors $T^i$ we refer to the papers cited above. In order to simplify notations we adopt to the conventions
\[T^i(B/A):=T^i(B/A,B),\qquad T^i(B,M):=T^i(B/\C,M).\]
Additionally, if $X\to S$ is a morphism of space germs and $M$ is an $\OO_X$-module we also write
\[T^i(X/S,M):=T^i(\OO_X/\OO_S,M).\]
\end{remark}

\begin{proposition}\label{tiprop}
The functors $T^i$ have the following properties:
\begin{enumerate}
\item $T^0(B/A,M)\simeq \Der_A(B,M)$.
\item If $A\to B$ turns $B$ into a regular $A$-algebra, then $T^i(B/A,M)=0$ for $i>0$.
\item If $B\simeq A\{\x\}/(f_1,\dots,f_k)$, where $f_1,\dots,f_k$ form a regular sequence in $A\{\x\}$, then $T^i(B/A,M)=0$ for $i>1$.
\item If $A\onto B$ is a surjection with kernel $I$ then $T^1(B/A,M)\simeq \Hom_B(I/I^2,M)$.
\item If $T^1(B/A,\C)=0$, then $B$ is a regular $A$-algebra.
\item If $T^2(B/A,\C)=0$, then $B\simeq A\{\x\}/(f_1,\dots,f_k)$ with a regular sequence $f_1,\dots,f_k\in A\{\x\}$.
\end{enumerate}
\end{proposition}

The functorial properties of cotangent cohomology yield the following exact sequence:

\begin{proposition}\label{exmodseq}
Let $0\to M'\to M\to M''\to 0$ be a short exact sequence of $B$-modules. Then there is a long exact sequence in cohomology
\[
0\to T^0(B/A,M')\to T^0(B/A,M)\to T^0(B/A,M'')\to T^1(B/A,M')\to \dots
\]
\end{proposition}

Furthermore, suppose we are additionally given a morphism $B\to C$, then there is an exact sequence
\[0\to \LL_{B/A}^\bullet\otimes_B C\to\LL_{C/A}^\bullet\to \LL_{C/B}^\bullet \to 0,\]
from which one derives:

\begin{proposition}\label{threealg}
If $A\to B\to C$ are homomorphisms of analytic algebras and $M$ is a $C$-module, then there is a long exact sequence in cohomology
\[
0\to T^0(C/B,M)\to T^0(C/A,M)\to T^0(B/A,M)\to T^1(C/B,M)\to \dots
\]
\end{proposition}

In particular, combining Proposition \ref{tiprop} and Proposition \ref{threealg} we get:

\begin{proposition}\label{normalsequence}
Let $B\simeq A\{\x\}/I$ be a quotient of a regular analytic $A$-algebra $A\{\x\}$, and let $M$ be a $B$-module. Then there is an exact sequence
\[0\to \Der_A(B,M)\to \Der_A(A\{\x\},M)\to \Hom_B(I/I^2,M)\to T^1(B/A,M)\to 0.\]
\end{proposition}

\begin{proof}
Applying Proposition \ref{threealg} to the homomorphisms $A\to A\{\x\}\to B$ we get the exact sequence
\[0\to T^0(B/A,M)\to T^0(A\{\x\}/A,M)\to T^1(B/A\{\x\},M)\to T^1(B/A,M)\to 0,\]
using property (2) in Proposition \ref{tiprop} and the fact that $T^0(B/A\{\x\},M)\simeq \Der_{A\{\x\}}(B,M)$ $=0$. Now the statement follows from the identities (1) and (4) in Proposition \ref{tiprop}.
\end{proof}

\begin{remark}
In fact, the sequence above results from applying $\Hom_B(-,M)$ to the so-called {\em conormal sequence} 
\[I/I^2\to \Omega_{A\{\x\}/A}\otimes_{A\{\x\}} B\to\Omega_{B/A}\to 0,\]
cf. \cite{Eisenbud}, for instance.
\end{remark}

Finally, the following proposition describes the behaviour of the $T^i$-functors under base change:

\begin{proposition}\label{basechange}
Let be given homomorphisms $A\to B$ and $A\to A'$ of analytic algebras, and set $B':=B\otimes_A A'$, i. e. we have a cocartesian diagram
\[
\xymatrix{B\ar[r]& B'=B\otimes_A A'\\
A\ar[u]\ar[r]&A'.\ar[u]}\]
If $A\to A'$ or $A\to B$ is flat, then there is a natural isomorphism
\[T^i(B/A,M)\simeq T^i(B'/A',M)\]
for any $B'$-module $M$ and $i\ge 0$.
\end{proposition}

\begin{remark}
If, in the situation above, both $A\to A'$ and $B\to B'$ are flat, then there is also the identity $T^i(B/A,M)\otimes_A A'\simeq T^i(B'/A',M\otimes_A A')$ for any $B$-module $M$. If only $A\to A'$ is flat, we still have such a map, but in general this map needs not be an isomorphism.
\end{remark}

So far we have introduced the notions and basic properties of the relative cotangent cohomology of a morphism $A\to B$ of analytic algebras; taking $A:=\C$ we get the absolute cotangent cohomology of $B$. Now let $X_0$ be a germ of a complex space, then one has (cf. \cite[§ 5]{Pal76}, \cite[2.4.4]{Buch81}) the following basic interpretation of $T^i(X_0)$ in terms of {\em deformations} of $X_0$, i. e. flat maps $\xi:X\to S$ with $\xi^{-1}(0)\simeq X_0$ (here and subsequently we denote by $D$ the one-point space $\{0\}$ with $\OO_D\simeq \C[\eps]/(\eps^2)$):
\begin{itemize}
\item $T^0(X_0)$ is isomorphic to the module of infinitesimal automorphisms of $X_0$, i. e. automorphisms of $X_0\times D$ over $D$,
\item $T^1(X_0)$ coincides with the space of first-order infinitesimal deformations, i. e. the tangent space of the {\em deformation functor}
\[Def_{X_0}(S):=\{\text{isomorphism classes of deformations of }X_0\text{ over }S\},\]
\item $T^2(X_0)$ contains the obstructions of lifting deformations to higher order. 
\end{itemize}
Analogous statements hold for the modules $T^i(X/S)$ that control relative deformations of $X$ over $S$, given by a diagram
\[\xymatrix{
X\ar[dd]\ar@{^{(}->}[r]^i &\mc{X}\ar[dd]\ar[dr]^{\xi} \\
&& B.\\
S\ar@{^{(}->}[r]^(0.4)j &S\times B\ar[ur]_{pr_2}}\]

\subsection{The cohomology functors $T^i(f,-)$}

In what follows we will also consider deformations of morphisms $f:X\to S$ of germs, their corresponding cotangent cohomology is given by functors $T^i(f,-)$ that we are now going to describe. Again, for details and the application of the following constructions (and their global counterparts) to deformations of morphisms of complex spaces, we refer to \cite{Flenn78}. 

\begin{definition}
Let $A,B$ and $C$ be analytic algebras such that $B$ and $C$ are $A$-algebras, and let $f:B\to C$ be an $A$-algebra homomorphism (for the applications we have in mind, we are of course mainly interested in the case $A=\C$). Choose resolvents $R$ (resp. $S$) of $B$ (resp. $C$) over $A$ that are compatible with $f$, i. e. there is a commutative diagram
\[\xymatrix{
R\ar[r]\ar[d]& S\ar[d]\\
B\ar[r]^f&C}\]
of DG-algebras over $A$. 
\begin{enumerate}
\item Define
\[\LL_f^\bullet:=(\LL_{B/A}^\bullet,\LL_{C/A}^\bullet).\]
Formally, this is a complex in the following category $\mc{C}$: $\mc{C}$ is the (abelian) category whose objects are tripels $(M_1,M_2,h)$ of a $B$-module $M_1$, a $C$-module $M_2$ and a homomorphism $h:M_1\to M_2$ over $f:B\to C$. $Mor_\mc{C}((M_1,M_2,h),(M'_1,M'_2,h'))$ consists of all pairs $(g_1,g_2)$ of a $B$-module homomorphism $g_1:M_1\to M'_1$ and a $C$-module homomorphism $g_2:M_2\to M'_2$ such that $g_2\circ h=h'\circ g_1$. $\LL_f^\bullet$ is then a complex in the derived category $D^-(\mc{C})$.
\item For a $B$-module $M$ we set
\[T^i(f/A,M):=\Ext_\mc{C}^i(\LL_f^\bullet,(M,M\otimes_B C)),\qquad i\ge 0.\]
As before, we simply write $T^i(f,M)$ if $A=\C$ and $T^i(f/A):=T^i(f/A,B)$.
\item More generally, we write
\[T^i(f/A,(M_1,M_2)):=\Ext^i_\mc{C}(\LL_f^\bullet,(M_1,M_2)),\qquad i\ge 0,\]
for any object $(M_1,M_2,h)$ in $\mc{C}$.
\end{enumerate}
\end{definition}

Again, one can describe $T^0(f,-)$ more explicitly in terms of derivations:

\begin{proposition}
Let $f:B\to C$ be as above and let $M$ be a $B$-module. Then $T^0(f/A,M)$ consists of all pairs of compatible derivations $(\delta,\tilde\delta)\in \Der_A(B,M)\times \Der_A(C,M\otimes_B C)$, i. e. pairs $(\delta,\tilde\delta)$ such that
\[\xymatrix{
C\ar[r]^{\tilde\delta}&M\otimes_B C\\
B\ar[r]^{\delta}\ar[u]^f&M\ar[u]
}\]
commutes.
\end{proposition}

Furthermore, there is, for any $B$-module $M$, a short exact sequence of complexes
\[0\to \Hom_C(\Omega_{S/R}\otimes_S C,M\otimes_B C)\to \Hom_\mc{C}(\LL_f^\bullet(M,M\otimes_B C))\to \Hom_B(\Omega_{R/A}\otimes_R B,M)\to 0\]
which yields the following long exact cohomology sequence:

\begin{proposition}\label{T-rel-morph}
For any morphism $f:B\to C$ of analytic $A$-algebras and any $B$-module $C$ there is a long exact sequence in cohomology
\[
0\to T^0(C/B,M\otimes_B C)\to T^0(f/A,M)\to T^0(B/A,M)\to T^1(C/B,M\otimes_B C)\to\dots
\]
\end{proposition}

\begin{corollary}\label{ks-sequence}
Let $\xi:X\to S$ be a morphism of complex space germs and $M$ an $\mc{O}_S$-module. Then:
\begin{enumerate}
\item There is a long exact sequence
\[0\to T^0(X/S,M\otimes_{\OO_S}\OO_X)\to T^0(\xi,M)\to T^0(S,M)\to 
T^1(X/S,M\otimes_{\OO_S}\OO_X)\to\dots\]
\item Taking $M:=\OO_S$ one obtains the so-called {\em Kodaira-Spencer-sequence} for the mapping $\xi$:
\[0\to T^0(X/S)\to T^0(\xi)\to T^0(S)\to T^1(X/S)\to\dots\]
In particular, for $i=0$ we have
\[T^0(\xi)=\left\{(\delta,\tilde\delta)\in \Der_\C(\OO_S,\OO_S)\times \ Der_\C(\OO_X,\OO_X) : \xi^*\circ\delta=\tilde\delta\circ \xi^*\right\},\]
where $\xi^*$ denotes the corresponding homomorphism $\OO_S\to\OO_X$.
\end{enumerate}
\end{corollary}

Again, for small $i$ the modules $T^if)$ have an interpretation in the context of deformations of a morphism $f:X\to S$: $T^0(f)$ contains the pairs of compatible infinitesimal automorphisms of $X$ and $S$ over the double point $D$, and $T^1(f)$ coincides with the set of isomorphism classes of deformations of $f$ over $D$, given by a commutative diagram
\[\xymatrix{
X\ar[dd]_{f}\ar@{^{(}->}[r]^i &\mc{X}\ar[dd]^{F}\ar[dr]^{\xi} \\
&& D,\\
S\ar@{^{(}->}[r]^j &\mc{S}\ar[ur]_{\rho}}\]
where $\xi$ and $\rho$ are deformations of $X$ and $S$, respectively.

%-- Deformations with section -------------------------------------------------

\section{Deformations with section}

In this section we review the definition of the functor $Def^s_{X_0}$  of deformations with  section of an isolated singularity $X_0$ and collect some results on the corresponding cotangent cohomology modules and the construction of (mini-)versal deformations with section, see also \cite{Buch81}.

\subsection{The functor $Def^s_{X_0}$}

\begin{definition}
Let $X_0$ be a germ of a complex space.
\begin{enumerate}
\item A {\em deformation with section} of $X_0$ over a space germ $S$ is a deformation $\xi:X\to S$ of $X_0$ together with a section $\sigma:S\to X$, i. e. $\xi\circ\sigma=id_{S}$. That is, we have a cartesian diagram
\[\xymatrix{
0\ar@{^{(}->}[r]\ar[d]&S\ar[d]^\sigma\\
X_0\ar@{^{(}->}[r]^i\ar[d]&X\ar[d]^\xi\\
0\ar@{^{(}->}[r]&S,
}\]
and the corresponding homomorphism $\xi^*:\mc{O}_{S}\to\mc{O}_{X}$ turns $\mc{O}_{X}$ into a flat $\mc{O}_{S}$-module. We write $(\xi,\sigma)$ for a deformation together with its section.\\
\item Two deformations $\xi:X\to S$ and $\xi':X'\to S$ over $S$ with sections $\sigma$ and $\sigma'$ are {\em isomorphic} if there is a morphism $\psi:X\to X'$ preserving $X_0$ inside $X$ resp. $X'$ and making the diagram
\[\xymatrix{
&S\ar[dl]_\sigma\ar[dr]^{\sigma'}\\
X\ar[dr]_\xi\ar[rr]^\psi&&X'\ar[dl]^{\xi'}\\
&S
}\]
commute (and $\psi$ is then an isomorphism, whence the terminology).
\item {\em Base change:} If $\ph:T\to S$ is a further map of space germs, the {\em pull-back} (or {\em induced deformation}) of $\xi$ by $\ph$ is defined as the fibre product of $X$ and $T$ over $S$ together with the obvious maps:
\[\xymatrix{
T\ar[r]^\ph\ar[d]_{t\mapsto (\sigma(\ph(t)),t)}&S\ar[d]^\sigma\\
X\ar[r]^(.6){pr_1}\ar[d]_{pr_2}\times_S T&X\ar[d]^\xi\\
T\ar[r]^\ph&S.
}\]
We denote the pull-back as $\ph^*(\xi)$; if $T\stackrel{i}{\hookrightarrow}S$ is a subspace we simply write $\xi_{|T}$ for $i^*(\xi)$.
\item With these notions, we can define the {\em functor of deformations with section of $X_0$}
\begin{align*}
Def_{X_0}^s:(Gan)&\to (Sets)\\
S&\mapsto
\left\{\BMA
\text{isomorphism classes of deformations}\\
\text{with section of }X_0\text{ over }S
\BME\right\}.
\end{align*}
By abuse of notation, we denote by $Def_{X_0}^s$ the (covariant) functor $(Analg)\to (Sets)$ of deformations with section of the analytic algebra $\mc{O}_{X_0}$, too.
\end{enumerate}
\end{definition}

\subsection{Cotangent cohomology of $Def^s_{X_0}$}

We can interpret deformations with section of $X_0$ as deformations of the morphism $0\hookrightarrow X_0$ corresponding to the residue map $\mc{O}_{X_0}\onto \C$. Its cotangent cohomology coincides with the cotangent cohomology of $X_0$ with values in $\mf{m}$, the maximal ideal of $\mc{O}_{X_0}$:

\begin{proposition}\label{t12-sec}
Let $X_0$ be a germ of a complex space, minimally embedded in $(\C^n,0)$ and given by an ideal $I=(f_1,\dots,f_k)\subseteq \OO_n$. Then  \[T^i(0\hookrightarrow X_0)=T^i(X_0,\mf{m}), i\ge 0,\]
where  $\mf{m}\subseteq \mc{O}_{X_0}$ is the maximal ideal. In particular we have:
\begin{enumerate}
\item $T^0(X_0,\mf{m})\simeq \Der_\C(\OO_{X_0},\mf{m})$.
\item $T^1(X_0,\mf{m})$ is the cokernel of the map
\begin{align*}
\Der_\C(\OO_n,\mf{m})&\to \Hom_{\OO_{X_0}}(I/I^2,\mf{m})\\
\delta&\mapsto (g+I^2\mapsto \delta(g)).
\end{align*}
If $X_0$ is a complete intersection and $f_1,\dots,f_k$ form a regular sequence, then
\[T^1(X_0,\mf{m})\simeq \mf{m}^{\oplus k}/(\mf{m}J(f)),\]
where $J(f)$ is the module generated by the columns of the Jacobian matrix $\left(\SMA
\overline{\frac{\p f_1}{x_1}} & \cdots & \overline{\frac{\p f_1}{x_n}}\\
\vdots&&\vdots\\
\overline{\frac{\p f_k}{x_1}} & \cdots & \overline{\frac{\p f_k}{x_n}}
\SME\right)$.
\item Complete intersections are unobstructed for $Def^s_{X_0}$, i. e. $T^2(X_0,\mf{m})=0$.
\end{enumerate}
\end{proposition}       
       
\begin{proof}
Let $r$ denote the residue map $\OO_{X_0}\onto\C$. The cotangent complex of $\C$ being trivial we get $\Hom_\mc{C}(\LL_r^\bullet,(\OO_{X_0},\C))\simeq \Hom_{\OO_{X_0}}(\LL_{\OO_{X_0/\C}}^\bullet,\mf{m})$, where $\mc{C}$ denotes the category described in section 1.2. Thus
\begin{align*}
T^i(0\hookrightarrow X_0)
&=\Ext^i_\mc{C}(\LL_r^\bullet,\C)\\
&\simeq \Ext_{\OO_{X_0}}^i(\LL_{\OO_{X_0/\C}}^\bullet,\mf{m})
=T^i(X_0,\mf{m}).
\end{align*}
(1) and (3) are a direct consequence from Proposition \ref{tiprop}. Property (2) follows from Proposition \ref{normalsequence}, noting that $\Der_\C(\OO_n,\mf{m})$ is generated by $\left\{x_i\frac{\p}{\p x_j}\right\}_{i,j}$ and $\Hom_{\OO_{X_0}}((f)/(f^2),\mf{m})\simeq \mf{m}^{\oplus k}$ for complete intersections. 
\end{proof}

The same argument as in the case of the functor $Def_{X_0}$ shows:

\begin{lemma}\label{infi-defs}
Let $X_0$ be a germ of a complex space. Then $T^1(X_0,\mf{m})$ is isomorphic to the $\OO_{X_0}$-module of isomorphism classes of first-order infinitesimal deformations with section of $X_0$, where $\mf{m}$ denotes the maximal ideal of $\OO_{X_0}$, i. e. $T^1(X_0,\mf{m})$ coincides with the tangent space of the functor $Def^s_{X_0}$.
\end{lemma}

\begin{remark}
In particular, as embedded deformation, every infinitesimal deformation of $X_0\subseteq (\C^n,0)$ given by $f_1,\dots,f_k\in\OO_n$ is defined by $f_i+\eps g_i\in \OO_n\otimes \C[\eps]/(\eps^2)$, $i=1,\dots,k$, where the $g_i$ are elements of the maximal ideal of $\OO_n$. More generally, by a suitable coordinate change we may assume that (up to isomorphism) every deformation with section of $X_0$ over the base $S$ is given by $F_1,\dots,F_k\in (\x)\OO_S\{x\}$, and $\sigma$ is the zero section $g\mapsto g\bmod{(\x)}$.
\end{remark}

\subsection{Construction of a versal deformation with section}

We now turn to the construction of a versal deformation with section of an isolated singularity $X_0$, which goes back to \cite{Buch81}. There are the usual notions of (mini-)versality:

\begin{definition}
Denote by $h_Y$ the $(Sets)$-valued functor $X\mapsto \Hom(X,Y)$.
\begin{enumerate}
\item A deformation $\xi:X\to S$ with section $\sigma$ is called {\em versal} (as deformation with section) if for any complex space germ $T$ the map
\[h_S(T)\to Def_{X_0}^s(T),\qquad \ph\mapsto \ph^*(\xi)\]
is surjective, i. e. every $(\nu,\tau)\in Def_{X_0}^s(T)$ is induced from $\xi$ by a morphism $T\to S$.
\item As usual, the natural map $\theta_\xi:T(S)\simeq \Hom(\OO_S,\C[\eps])\to Def_{X_0}^s(D)$ between the tangent spaces of the functors $h_S$ and $Def_{X_0}^s$ is called the {\em Kodaira-Spencer map} of the deformation. $(\xi,\sigma)$ is called {\em semi-universal} or {\em miniversal} (as deformation with section), if it is versal and, in addition, $\theta_\xi$ is a bijection.
\end{enumerate}
\end{definition}

\begin{theorem}\label{versal-sec}
Suppose $\xi:X\to S$ is a versal deformation of $X_0$. Then $\pi_1:X\times_S X\to X$ together with the diagonal embeding $d:X\to X\times_S X$ gives a versal deformation with section of $X_0$.
\end{theorem}

\begin{proof}[\cite{Buch81}, see also \cite{MvS01}] Let $\nu:\tilde X\to T$ be a deformation of $X_0$ with section $\tau$. Versality of $\xi$ implies the existence of some $\ph:T\to S$ such that $\nu$ is isomorphic to the pull-back of $\xi$ via $\ph$, so we may assume $\tilde X=X\times_S T$. Denoting by $\Phi$ the projection $\tilde X\to X$ we obtain by pulling back $X\times_S X\to X$ over $\Phi\circ \tau$:
\[\xymatrix{
& X\times_S X\ar[r]^{\pi_2}\ar[dd]^{\pi_1}&X\ar[dd]^\xi\\
(X\times_S X)\times_X T\ar[ur]\ar[dd]\\
&X\ar[r]^\xi&S\\
T\ar[ur]^{\Phi\circ\tau}\ar[r]^\tau&\tilde X\ar[r]^\nu\ar[u]_\Phi&T\ar[u]_\ph
}\]
From $\tilde X=X\times_S T$ we deduce that (up to isomorphism):
\[(X\times_S X)\times _X T\simeq X\times_S(X\times_X T)\simeq X\times_S T=\tilde X.\]
Under these identifications, we can identify $d\circ (\Phi\circ\tau)$ with $\tau$. Altogether this proves that $(\nu,\tau)$ is isomorphic to the pull-back of $(\pi_1,d)$ via $\Phi\circ\tau$.
\end{proof}

Since, by the classical result in \cite{Gra72}, any isolated singularity admits a versal deformation, we obtain immediately:

\begin{corollary}
Every isolated singularity $X_0$ has a (mini-)versal deformation with section.
\end{corollary}

\begin{example}[Complete intersections]\label{sud-sec-icis}
We give an explicit construction of a versal deformation of an isolated complete intersection singularity (ICIS) in terms of a defining regular sequence $f=(f_1,\dots,f_k)\in\C\{x_1,\dots,x_n\}^k$: Take a family $\{g^{(1)},\ldots,g^{(\tau)}\}$ of monomials (i. e. $g^{(i)}=x^{\alpha_i}e_{k(i)}$, where $e_j$ denotes the $j$-th unit vector) representing a $\C$-vector space basis of $T^1(X_0)$, we may assume that $g^{(i)}=e_i$ for $i=1,\ldots,k$. Then
\begin{multline*}
\OO_S:=\C\{u,s\}/(F_1(u,s),\dots,F_k(u,s))  \\
\to
\C\{x,u,s\}/
(F_1(x,s),\dots,F_k(x,s),F_1(u,s),\dots,F_k(u,s))=:\OO_X
\end{multline*}
is a versal deformation with section of $X_0$, where
\[\left(\SMA F_1(x,s)\\\vdots\\F_k(x,s)\SME\right):=
\left(\SMA f_1(x)\\\vdots\\f_k(x)\SME\right)+
\sum_{i=1}^\tau s_ig^{(i)}(x)\]
and the section is given by
\[\OO_X\to \OO_S,\quad s_i\mapsto s_i,\ u_i\mapsto u_i,\ x_i\mapsto u_i.\]
$\OO_S$ is regular, using the relation
\[\left(\SMA s_1\\\vdots\\s_k\SME\right)=-\left(\SMA f_1(u)\\\vdots\\f_k(u)\SME\right)-\sum_{i=k+1}^\tau s_ig^{(i)}(u)\]
in $\OO_X$ one obtains that this is isomorphic to
\[\C\{u,s_{k+1},\dots,s_{\tau}\}\to
\C\{x,u,s_{k+1},\dots,s_{\tau}\}/
(\tilde F_1(x,u,s),\dots,\tilde F_k(x,u,s)),\]
where
\[\left(\SMA \tilde F_1(x,u,s)\\\vdots\\ \tilde F_k(x,u,s)\SME\right):=
\left(\SMA f_1(x)-f_1(u)\\\vdots\\f_k(x)-f_k(u)\SME\right)
+\sum_{i=k+1}^\tau s_i(g^{(i)}(x)-g^{(i)}(u)).\]
From the following lemma we can conclude that this deformation is also semi-universal. 

At the end of this section a different construction of a semi-universal deformation using a $\C$-basis of $T^1(X_0,\mf{m})$ is given.
\end{example}

\begin{lemma}\label{dimT1sec}
Let $f=(f_1,\dots,f_k)\in\OO_n^k$ define an ICIS $X_0:=(V(f),0)$. Then
\[\tau^s(X_0):=\dim_\C T^1(X_0,\mf{m})=\tau(X_0)+n-k,\]
where $\tau(X_0)=:\dim_\C T^1(X_0)$ denotes the Tjurina number of $X_0$.
\end{lemma}

\begin{proof}
Denote by $(\x)$ the maximal ideal of $\OO_n=\C\{x_1,\dots,x_n\}$. By Lemma \ref{t12-sec} we have \[T^1(X_0,\mf{m})=\mf{m}^{\oplus k}/(\mf{m}J(f))=(\x)^{\oplus k}/((f_i e_j)_{i,j}+(x_i\tfrac{\p f}{\p x_j})_{i,j}).\]
Since $(\x)^{\oplus k}$ has $\C$-codimension $k$ in $\OO_n^k$, the above formula is equivalent to the $\C$-linear independency of the $n$ partial derivatives $\frac{\p f}{\p x_i}$ of $f$ modulo $(f,(\x)J(f))$. Suppose the converse, i. e. there exist $\lambda_1,\dots,\lambda_n\in\C$, not all $\lambda_i$ equal to zero, such that $\sum_i \lambda_i \frac{\p f}{\p x_i}\equiv 0\ ((f_i e_j)_{i,j}+(\x)J(f))$. After a linear change of coordinates we may assume that $\frac{\p f}{\p x_1}\equiv 0\ ((f_i e_j)_{i,j}+(\x)J(f))$, i. e.
\[\tfrac{\p f_i}{\p x_1}=bf_i+\sum_{j=1}^n c_j\tfrac{\p f_i}{\p x_j},\quad b\in\OO_n,\ c_j\in (\x)\qquad (*)\]
for all $i$. In addition, by a suitable choice of a generating system of the ideal $(f_1,\dots,f_k)$, we may assume that each generator $f_i$ defines an isolated singularity, hence without loss of generality we may also assume that they are of the form
\[f_i= a_{i1}x_1^{d_{i1}}+\sum_{j=2}^n a_{ij}x_1^{d_{ij}}x_j+r_i(x),\quad a_{ij}\in\C,\ r_i\in (x_2,\dots,x_n)^2,\ i=1,\dots,n,\]
and not all $a_{ij}=0$, $j=1,\dots,n$. But taking the partial derivative with respect to $x_1$ then immediately gives a contradiction to $(*)$.
\end{proof}

For arbitrary isolated singularities $X_0$ defined by $f_1,\dots,f_k\in \OO_n$ one may take a different approach to determine (at least $k$-jets of) a (mini)-versal deformation with section: First compute $g_1,\dots,g_{\tau^s}\in\OO_n^k$ representing a $\C$-vector space basis of $T^1(X_0,\mf{m})$. We obtain a versal family of first order deformations
\[F^{(1)}:=f+\sum_{i=1}^{\tau^s} s_i g_i\in\left(\OO_n\otimes \C\{s\}/(s)^2\right)^k\]
over $\C\{s\}/(s)^2$. Then one lifts this family order by order (in $s$), killing obstructions. This is the algorithm described in \cite{Mar98} for the computation of versal deformations of singularities, its implementation can be found in the \Singular-library \verb+deform.lib+ (\cite{deform}). 

Since complete intersection are unobstructed, we obtain that in particular a semi-universal deformation with section of an ICIS $X_0$ is given by
\[F:=f+\sum_{i=1}^{\tau^s}s_i g_i\in \OO_{n+\tau^s}^k.\]

%-- Modular deformations ------------------------------------------------------

\section{Modular deformations}

\subsection{Definition and basic properties}

A deformation (with section) $\xi: X\to S$ of an isolated singularity $X_0$ is called {\em universal} if every deformation $X'\to T$ of $X_0$ is (up to isomorphism) induced by a unique morphism $T\to S$, i. e. $\xi$ is versal and any two different morphisms $\ph,\psi:T\to S$ induce non-isomorphic deformations of $X_0$. While any isolated singularity admits a semi-universal deformation (with section), the occurence of trivial subfamilies in such deformations implies that, in general, universal deformations of $X_0$ cannot exist, i. e. the functors $Def_{X_0}$ and $Def^s_{X_0}$ are not representable. This is, for instance, the case if $X_0$ is an isolated complete intersection singularity.

 However, restricting a semi-universal familiy to subgerms of the base for which the universality condition holds is a possible approach to the construction of local moduli for singularities. We call such subgerms {\em modular}, as introduced by Palamodov in \cite{Pal78}. The similar notion of {\em prorepresenting substratum} inside the base space of a semi-universal deformation of $X_0$ is -- in the context of formal deformations -- studied in \cite{LP88}.

\begin{definition}
Let $\xi:X\to S$ be a deformation of $X_0$ (with section $\sigma$). A subspace $M\subseteq S$ is called {\em modular} if the following condition holds: If $\ph:T\to M$ and $\psi:T\to S$ are morphisms such that the induced deformations $\ph^*(\xi_{|M})$ and $\psi^*(\xi)$ with base $T$ are isomorphic as deformations (with section) then $\ph=\psi$.

The restriction of a deformation to a modular subgerm is called a {\em modular deformation}.
\end{definition}

\begin{remark}
The following properties are immediate from the definition:
\begin{enumerate}
\item Any subgerm of a modular germ is again modular.
\item $\{0\}\subseteq S$ is a modular subspace for the deformation $\xi:X\to S$ (with section $\sigma$) if and only if the corresponding Kodaira-Spencer map $\theta_\xi: T(S)\to T^1(X_0)$ (resp. $\theta_\xi: T(S)\to T^1(X_0,\mf{m})$) is injective. Together with (1) this implies that $S$ contains a modular subspace if and only if the Kodaira-Spencer map of the deformation is injective. Following \cite{Pal90} we call such deformations {\em monodeformations}.
\item If $M_1,M_2$ are modular subgerms of $S$, then $M_1\cup M_2\subseteq S$ is modular, too.
\item As a consequence of the identity theorem for power series, two morphisms between space germs $\ph,\psi:T\to S$ coincide if $\ph_{|W}=\psi_{|W}$ for any Artinian subgerm $W$ of $T$, thus it suffices to check the condition of the definition for Artinian germs $T$.
\item There is a unique isomorphism between any two maximal modular subgerms in the base space of a semi-universal deformation of $X_0$. Such a maximal  modular subgerm is then called the {\em modular stratum} of $X_0$.
\end{enumerate}
\end{remark}

\subsection{The Kodaira-Spencer sequence of a deformation}

In order to derive a criterion for $M\subseteq S$ to be modular, we consider the following commutative diagram with exact rows, combining the Kodaira-Spencer sequence for $\xi:X\to S$ (Corollary \ref{ks-sequence} (2)) and its evaluation at the special fibre over $0\in S$:
\begin{gather}\label{kss}\xymatrix{
0\ar[r] & T^0(X/S)\ar[r]\ar[d]^{ev} & T^0(\xi)\ar[r]\ar[d]^{ev'} & T^0(S)\ar[r]^(0.45){\Theta_\xi}\ar[d]
& T^1(X/S)\ar[r]\ar[d] & \dots \\
0\ar[r] & T^0(X_0)\ar[r] & T^0(\xi,(\C,\OO_{X_0}))\ar[r] & T^0(S,\C)\ar[r]^{\theta_\xi}
& T^1(X_0)\ar[r] & \dots,
}\end{gather}
cf. \cite[§ 1]{Pal90}. Evaluation at the special fibre means taking cohomology with values in
$\C\simeq\OO_S/\mf{m}_S$ instead of $\OO_S$. Furthermore in the lower row we apply Proposition \ref{basechange}.
to obtain $T^i(X/S,M)\simeq T^i(X_0,M)$ for any $i\ge 0$ and any $\OO_{X_0}$-module $M$. Here $\theta_\xi$ denotes the {\em Kodaira-Spencer map} of the deformation, and $\Theta_\xi$ is the corresponding {\em relative Kodaira-Spencer map} $T^0(S)\to T^1(X/S)$.\\

The following proposition gives an analogue of sequence (\ref{kss}) for deformations with section.

\begin{proposition}\label{ks-modif}
Let $\xi:X\to S$ be a deformation of $X_0$ with section $\sigma$, $J_\sigma\subseteq \OO_X$ the kernel of $\sigma^*$. Then the sequence
\begin{align}\label{ks-sequence-modif}\xymatrix{
0\ar[r]& T^0(X/S,J_\sigma)\ar[r]& T^0(\xi,(\OO_S,J_\sigma))\ar[r]& T^0(S)\ar[r]^(0.4){\Theta_\xi} &T^1(X/S,J_\sigma)\ar[r]&\ldots
}\end{align}
is exact and by evaluation at the special fibre over $0\in S$ we obtain the diagram with exact rows
\begin{align}\label{ks-diag-modif}\xymatrix{
0\ar[r]& T^0(X/S,J_\sigma)\ar[d]_{ev}\ar[r]& T^0(\xi,(\OO_S,J_\sigma))\ar[d]_{ev'}\ar[r]& T^0(S)\ar[d]\ar[r]^(0.4){\Theta_\xi} &T^1(X/S,J_\sigma)\ar[d]\ar[r]&\ldots\\
0\ar[r]& T^0(X_0,\mf{m})\ar[r]& T^0(\xi,(\C,\mf{m}))\ar[r]& T^0(S,\C)\ar[r]^(0.45){\theta_\xi} &T^1(X_0,\mf{m})\ar[r]&\ldots ,
}\end{align}
where $\mf{m}$ denotes the maximal ideal of $\OO_{X_0}$. We call (\ref{ks-sequence-modif}) the {\em Kodaira-Spencer sequence} and $\Theta_\xi$ the {\em relative Kodaira-Spencer map} of $(\xi,\sigma)$.
\end{proposition}

\begin{proof}
Choose resolvents $R_X$ and $R_S$ of $\OO_X$ resp. $\OO_S$ over $\C$ as in the definition of the cotangent complex $\LL^\bullet_{\xi^*}$ in section 1.2. The section $\sigma^*$ yields a split exact sequence
\[\xymatrix{0\ar[r] & J_\sigma \ar@<0.3ex>[r]^\iota & \OO_X\ar@<0.3ex>[l]^\tau\ar@<0.3ex>[r]^{\sigma^*} & \OO_S\ar@<0.3ex>[l]^{\xi^*}\ar[r] & 0,}\]
hence we obtain a homomorphism $\tau\circ\xi^*:\OO_S\to J_\sigma$ over $\xi^*$, i. e. $(\OO_S,J_\sigma,\tau\circ\xi^*)$ is an object in the category $\mc{C}_{\xi^*:\OO_S\to\OO_X}$ defined in section 1. All cotangent complexes fit together into the short exact sequence
\[0\to \Hom_{\OO_X}(\Omega_{R_X/R_S}\otimes_{R_X} \OO_X,J_\sigma)\to \Hom_{\mc{C}}(\LL^\bullet_{\xi^*},(\OO_S,J_\sigma,\tau\circ\xi^*))\to \Hom_{\OO_S}(\LL^\bullet_{\OO_S/\C},\OO_S)\to 0.\]
The choice of $R_X$ and $R_S$ implies that $\Omega_{R_X/R_S}\otimes_{R_X} \OO_X\simeq \LL^\bullet_{\OO_X/\OO_S}$ (cf. \cite[ch. I]{Flenn78}), hence we get (\ref{ks-sequence-modif}) as its long exact cohomology sequence.

Similarly, if one takes coefficients in $(\C,\mf{m})$, this produces the exact  lower row of (\ref{ks-diag-modif}), and the commutativity of this diagram is then clear from functoriality.
\end{proof}

\subsection{A criterion for modularity in terms of $T^0$}

We are now going to give a criterion for $M\subseteq S$ to be modular, similar statements for deformations of complex spaces are already contained in \cite{Pal78}. In \cite{KS90} modular subspaces for general deformation groupoids are studied, and an analogue for criterion (3) below is given, involving a so-called {\em exponential functor} that generalizes the module of derivations $T^0(X/S)$. We use the following two auxiliary notions, the first of which being motivated by the fact that $\{0\}\subseteq S$ is modular if and only if $\theta_\xi$ is injective.

\begin{definition}
A subgerm $M\subseteq S$ in the base space of a deformation $\xi:X\to S$ (with section $\sigma$) of $X_0$ is called
\begin{enumerate}
\item {\em infinitesimally modular} if the restriction to $M$ of the map
\[\tilde\Theta_{\xi|M}:T^0(S,\OO_M)\to T^1(X/S)_{|M}\]
(resp. $\tilde\Theta_{\xi|M}:T^0(S,\OO_M)\to T^1(X/S,J_\sigma)_{|M}$ if $\xi$ is considered as deformation with section $\sigma$) is injective;
\item {\em Artinian modular} if any Artinian subgerm of $M$ is modular.
\end{enumerate}
(Here and subsequently we use the notations $T^i(X/S)_{|M}:=T^i(X_{|M}/M)$ etc. for the cotangent cohomology of the restricted deformation $\xi_{|M}$.)
\end{definition}

\begin{lemma}\label{mod-infmod}
For both deformation functors, modular subspaces are infinitesimally modular.
\end{lemma}

\begin{proof}[cf. \cite{Mar02}]
Suppose $\delta \in Ker(\tilde\Theta_{\xi|M})\subseteq T^0(S,\OO_M)$. Then $\delta$ corresponds to a morphism $\ph_\delta: M\times D\to S$ such that the deformation $\ph_\delta^*(\xi)$ is isomorphic to the trivial one, which can also be induced from $\xi_{|M}$ via the projection $pr_M:M\times D\to M$. By the modularity of $M$ we conclude that $\ph_\delta=i\circ pr_M$, hence $\delta=0$.
\end{proof}

\begin{lemma}\label{evprime-surj}
Suppose $\xi:X\to S$ is either a semi-universal or a modular deformation of $X_0$ (with section $\sigma$). Then, in each case, the mapping $ev':T^0(\xi)\to T^0(\xi,(\C,\OO_{X_0}))$ (resp. $ev':T^0(\xi,(\OO_S,J_\sigma))\to T^0(\xi,(\C,\mf{m}))$ is surjective.
\end{lemma}

\begin{proof}
First suppose $\xi$ is semi-universal. Take a derivation $\delta\in T^0(\xi,(\C,\mf{m}))$. Since $\xi$ is semi-universal, we may interpret $\delta$ as element of $T^0(X_0,\mf{m})$ and it suffices to find a preimage in $T^0(\xi,(\OO_S,J_\sigma))$.

We use the construction of \cite[Proposition 18.9]{Pal90}: $\delta$ induces an automorphism $a_\delta$ of $X_0\times D$ over $D$. Denote by $\tilde \xi$ the deformation  $\xi\times id_D:X\times D\to S\times D$ with section $\sigma\times id_D$. Let $T:=(\{0\}\times D)\cup(S\times\{0\})$ and $\ph:T\to S$ be the canonical projection onto $S$, so $\tilde \xi_{|T}=\ph^*(\xi)$. Let $\tilde a_\delta$ be the automorphism of $\ph^*(\xi)$ induced by $a_\delta$ on the first component of $T$ and by $id_S$ on the second. Then $\xi':=\tilde a_\delta\circ \ph^*(\xi)$ is of course still a versal deformation with section (isomorphic to $\tilde \xi_{|T}$), hence $\tilde \xi$ is induced from it by some $\psi:S\times D\to X$, which altogether yields the diagram
\[\xymatrix{
a' : X\times D\ar[r]\ar[d]^{\tilde\xi}&(X\times D)_{|T} \ar[rr]^(0.6){pr_X\circ \tilde a_\delta}\ar[d]^{\xi'} && X\ar[d]^\xi\\
\psi': S\times D\ar[r]^(0.6)\psi & T\ar[rr]^\ph && S
}\]
The corresponding map $\psi'^*:\OO_S\to\OO_S[\eps]$ gives a vector field $\eta\in \Der_\C(\OO_S)$ defined by $\psi'^*(a)=a+\eps\eta(a)$, similarly $a'$ induces a vector field $\tilde \delta\in \Der_\C(\OO_X,J_\sigma)$, with the property that $(\tilde \delta,\eta)\in T^0(\xi,(\OO_S,J_\sigma))$ and $ev'(\tilde \delta,\eta)=\delta$.

Now let $\hat\xi:Y\to M$ be a modular deformation with section $\hat\sigma$. We may assume that $M$ is a subspace of the base of a semi-universal deformation with section $\xi:X\to S$ of $X_0$, and $\hat\xi=\xi_{|M}$. As before, let $\tilde\xi:=\hat\xi\times id_D$, then $\delta\in T^0(\hat\xi,(\C,\mf{m}))\simeq T^0(X_0,\mf{m})$ induces an automorphism $\tilde a_\delta$ of the (still versal) deformation $\tilde\xi_{|T}$, so we can induce $\tilde\xi:Y\times D\to M\times D$ from $\xi':=a_\delta\circ \tilde\xi_{|T}$ by means of some morphism $\psi$:
\[\xymatrix{
a' : Y\times D\ar[r]\ar[d]^{\tilde\xi}&(X\times D)_{|T} \ar[rr]^(0.6){pr_X\circ \tilde a_\delta}\ar[d]^{\xi'} && X\ar[d]^\xi\\
\psi': M\times D\ar[r]^(0.6)\psi & T\ar[rr]^{\ph=pr_S} && S
}\]
On the other hand we can obtain $\tilde\xi$ from $\xi$ by just taking the projection $pr_1:M\times D\to M\subseteq S$. Thus, by the modularity of $M$, $pr_1=\psi'$, so that $a'$ and $\psi'$ can be considered as maps $Y\times D\to Y$ resp. $M\times D\to M$, inducing the desired lift of $\delta$ to $\Der_\C(\OO_Y,J_{\hat\sigma})$ resp. $\Der_\C(\OO_M)$.
\end{proof}

With these preparations we can formulate and prove the following criterion for modularity:

\begin{theorem}\label{modcrit}
For a subspace $M\subseteq S$ of the base space of a monodeformation $\xi:X\to S$ of $X_0$ the following assertions are equivalent:
\begin{enumerate}
\item $M$ is modular.
\item $M$ is infinitesimally modular.
\item $ev_{|M}:T^0(X/S)_{|M}\to T^0(X_0)$ is surjective.
\item $M$ is Artinian modular.
\end{enumerate}
\end{theorem}

\begin{theorem}\label{modcritsec}
For a subspace $M\subseteq S$ of the base space of a monodeformation $\xi:X\to S$ of $X_0$ with section $\sigma:S\to X$ the following assertions are equivalent:
\begin{enumerate}
\item $M$ is modular.
\item $M$ is infinitesimally modular.
\item $ev_{|M}:T^0(X/S,J_\sigma)_{|M}\to T^0(X_0,\mf{m})$ is surjective.
\item $M$ is Artinian modular.
\end{enumerate}
\end{theorem}

We only give a proof of Theorem \ref{modcritsec} -- Theorem \ref{modcrit} follows by omitting any section ocurring and replacing any $T^i(X_0,\mf{m})$ by $T^i(X_0)$, $T^i(X/S,J_\sigma)$ by $T^i(X/S)$ etc.

\begin{proof}[of Theorem \ref{modcritsec}]
$(1)\Rightarrow (2):$ This is Lemma \ref{mod-infmod}.

$(2)\Rightarrow (3):$ Consider the Kodaira-Spencer sequence (\ref{ks-sequence-modif}) of $(\xi,\sigma)$ with coefficients in $(\OO_M,J_{\sigma|M})$. This sequence can be put in between the rows of (\ref{ks-diag-modif}) so that, altogether, we obtain the commutative diagram with exact rows
\begin{align}\label{kss-subgerm}
\xymatrix{
0\ar[r]& T^0(X/S,J_\sigma)\ar[d]_{\tilde ev}\ar[r]& T^0(\xi,(\OO_S,J_\sigma))\ar[d]_{\tilde ev'}\ar[r]& T^0(S)\ar[d]\ar[r]^(0.45){\Theta_\xi} &T^1(X/S,J_\sigma)\ar[d]\ar[r]&\ldots\\
0\ar[r]& T^0(X/S,J_\sigma)_{|M}\ar[d]_{ev_{|M}}\ar[r]& T^0(\xi,(\OO_M,J_{\sigma|M}))\ar[d]_{ev'_{|M}}\ar[r]& T^0(S,\OO_M)\ar[d]\ar[r]^(0.45){\tilde\Theta_{\xi|M}} &T^1(X/S,J_\sigma)_{|M}\ar[d]\ar[r]&\ldots\\
0\ar[r]& T^0(X_0,\mf{m})\ar[r]& T^0(\xi,(\C,\mf{m}))\ar[r]& T^0(S,\C)\ar[r]^{\theta_\xi} &T^1(X_0,\mf{m}),\ar[r]&\ldots
}\end{align}
By assumption $\tilde\Theta_{\xi|M}$ is injective, $\theta_{\xi|M}$ is injective since $\xi$ is a monodeformation. Thus in (\ref{kss-subgerm}) we can identify $ev_{|M}$ and $ev'_{|M}$, the latter being surjective since, by Lemma \ref{evprime-surj}, the composition $ev'=ev'_{|M}\circ \tilde ev'$ is surjective.

$(3)\Rightarrow (4):$ Suppose $ev_{|M}$ is surjective and take an Artinian subgerm $A\subseteq M$. Since $A$ is obtained from $\C$ by a finite number of small extensions and $\{0\}\subseteq S$ is modular if and only if $\theta_\xi$ is injective, it suffices to prove the following statement: Let $\xi$ be modular over the Artinian subgerm $A_0$ of $M$, then it is also modular over $A$ given by an infinitesimal extension
\[\eta:0\to K\to \OO_A\to \OO_{A_0}\to 0.\]
As in \cite[Proposition 1.4]{Mar03} one associates to $\eta$ the obstruction element
\[ob_{\xi,\eta}\in \Hom(T^0(X_0,\mf{m}),T^1(X_0,\mf{m}))\otimes K,\]
whose vanishing is then equivalent to $A$ being modular, which is proved absolutely analogously to \cite[Lemma 1.5]{Mar03}. In concrete terms, we may assume that $\OO_A=\OO_{A_0}[\eps]$ and $K=(\eps)$, $\eps^2=0$. Furthermore let $X_0$ be given by $f=(f_1,\dots,f_k)\in \OO_n^k$, $X$ be defined by $F=(F_1,\dots,F_k)\in (\OO_n\otimes \OO_S)^k$.

Then we can write $F_{|A}=F_{|A_0}+\eps \tilde f$ with $\eps\tilde f\in \OO_n^k$ defining a class in $T^1(X_0,\mf{m})\otimes K$. Any $\delta\in T^0(X_0,\mf{m})$ is represented by $\delta'\in Der_\C(\OO_n,(x))$ with $\delta'(f)=h\cdot f$ for some matrix $h$ with entries in $\OO_n$. By assumption $(3)$ we can lift $\delta$ to some $\delta_{|A_0}\in T^0(X_{|A_0},J_{\sigma|A_0})$, represented by $\delta'_{|A_0}$ with $\delta'_{|A_0}(F_{|A_0})=H_{|A_0}F_{|A_0}$, $H_{|A_0}$ being a lift of $h$ to $\OO_n\otimes \OO_{A_0}$. A lift of $\delta'_{|A_0}$ and $H_{|A_0}$ to $\delta'_{|A}$ and $H_{|A}$ over $A$ induces classes
\[\delta'_{|A}(F_{|A})-H_{|A}F_{|A}=(\delta'(f)-h\tilde f)\eps\in T^1(X_0,\mf{m})\otimes K,\]
which yields a homomorphism $ob_{\xi,\eta}:T^0(X_0,\mf{m})\to T^1(X_0,\mf{m})\otimes K$. In fact, $ob$ is induced by the Lie bracket $[-,-]:T^0(X_0,\mf{m})\times T^1(X_0,\mf{m})\to T^1(X_0,\mf{m})$ which implies that it does not depend on the choices made.

Now by assumption (3) there exists a lift $\delta'_{|A}\in T^0(X_{|A},J_{\sigma|A})$, thus $ob_{\xi,\eta}$ vanishes, hence $A$ is modular.

$(4)\Rightarrow (1)$: Let $Z$ be an Artinian germ and $\ph:Z\to M$, $\psi:Z\to S$ morphisms such that $\ph^*(\xi_{|M})\simeq \psi^*(\xi)$. We can factor $\ph=i\circ\ph_0$, where $\ph_0:Z\to M_0$ with $M_0$ Artinian and $i:M_0\hookrightarrow M$ the inclusion. Thus $\psi^*(\xi)\simeq \ph^*_0(i^*(\xi_{|M}))\simeq\ph^*_0(\xi_{|M_0})$. By the assumption $M_0$ is modular, so $\psi=\ph_0$ and therefore $\psi$ and $\ph$ coincide, too.
\end{proof}

\begin{remark}
The definition of infinitesimal modularity used here differs slightly from the one used in \cite{Mar02}, \cite{Mar03}, where it is characterized as injectivity of the Kodaira-Spencer map $\Theta_{\xi|M}:T^0(M)\to T^1(X/S)_{|M}$ of the deformation $\xi_{|M}$. For any subspace $M\subseteq S$ we have a commutative diagram
\[\xymatrix{
T^0(S,\OO_M)\ar[rr]^{\tilde\Theta_{\xi|M}} && T^1(X_{|M}/M), \\
T^0(M)\ar[urr]_{\Theta_{\xi|M}}\ar@{^{(}->}[u]}\]
so $Ker(\tilde\Theta_{\xi|M})=0$ is the stronger notion, and it is indeed this condition that we need in the above proof of the implication $(2)\Rightarrow (3)$.
\end{remark}

As a corollary of the obstruction calculus in the proof of implication $(3)\Rightarrow (4)$ we obtain:

\begin{corollary}\label{tangmod}
The tangent space $T(M)$ to the modular stratum of $X_0$ equals the subspace
\[\{ t\in T(S):\ [\theta_\xi(t),\delta]=0\text{ for all } \delta\in T^0(X_0)\}\]
(resp. $T^0(X_0,\mf{m})$ for deformations with section).
\end{corollary}

\subsection{The module $T^0_\bullet(X_0)$}

Obviously, surjectivity of the evaluation mappings in the above statements is equivalent to the surjectivity of the mapping
\[T^0(X/S)_{|M}\to T^0(X_0)/ev(T^0(X/S))\quad \text{resp. } T^0(X/S,J_\sigma)_{|M}\to T^0(X_0,\mf{m})/ev(T^0(X/S,J_\sigma)).\]
This gives rise to the following notions (cf. \cite{Pal90}):

\begin{definition}\label{T0-red-def}
Let $\xi:X\to S$ be a semi-universal deformation of $X_0$ (with section $\sigma$). We set
\begin{itemize}
\item $T^0_\bullet(X_0):=T^0(X_0)/ev(T^0(X/S))$, resp.
\item $T^0_\bullet(X_0,\mf{m}):=T^0(X_0,\mf{m})/ev(T^0(X/S,J_\sigma))$ for the case of deformations with section.
\end{itemize}
\end{definition}

We note some properties of $T^0_\bullet(X_0)$ stated in \cite{Pal90}:
\begin{enumerate}
\item $\dim_\C T^0_\bullet (X_0)<\infty$;
\item $\dim_\C T^0_\bullet (X_0)\le \tau(X_0)$, with equality if and only if the kernel of the Kodaria-Spencer map $T^0(S)\to T^1(X/S)$ is free. This is, for instance, the case if $X_0$ is a complete intersection (by the Saito-Looijenga-Theorem, see \cite{Sai80}, \cite{Loo84}) or if $X_0$ is a reduced space curve singularity, as was shown by D. van Straten in \cite{Stra95};
\item if $X_0$ is a hypersurface defined by $f\in \OO_n$, then there is an exact sequence of $\OO_n$-modules
\begin{align}\label{T0red-seq}
0\to T^0_\bullet(X_0)\stackrel{\eps}{\to}\OO_n/J(f)\stackrel{m_f}{\to}
\OO_n/J(f)\stackrel{\pi}{\to}T^1(X_0)\to 0,
\end{align}
where $\pi$ is the projection, $m_f$ denotes multiplication by $f$ modulo $J(f)$ and $\eps$ is defined as follows: If $\delta\in T^0(X_0)$ is the class of a derivation $\delta':\OO_n\to\OO_n$ subject to $\delta'(f)=h\cdot f$ for some $h\in \OO_n$, then define $\eps(\delta):=\frac{\delta'(f)}{f}\bmod{J(f)}$. 
\end{enumerate}
 
Properties similar to (1) and (3) above hold for the module $T^0_\bullet(X_0,\mf{m})$. 
 
\begin{proposition}
Let $\xi:X\to S$ be a semi-universal deformation of $X_0$ with section $\sigma$. Then:
\begin{enumerate}
\item $\dim_\C T^0_\bullet(X_0,\mf{m})<\infty$;
\item if $X_0$ is a hypersurface defined by $f\in \OO_n$, then there is an exact sequence of $\OO_n$-modules
\begin{align}\label{T0red-sec-seq}
0\to T^0_\bullet(X_0,\mf{m})\stackrel{\eps}{\to}\OO_n/J(f)\stackrel{m_f}{\to}
\OO_n/(\x)J(f)\stackrel{\pi}{\to}\OO_n/(f,(\x)J(f))\to 0,
\end{align}
where the maps $\eps, m_f$ and $\pi$ are defined analogously as in (\ref{T0red-seq}). In particular we have:
\[\dim_\C T^0_\bullet(X_0,\mf{m})=\dim_\C T^0_\bullet(X_0)=\tau.\]
\end{enumerate}
\end{proposition}

\begin{proof}
Property (1) is a consequence of the Kodaira-Spencer sequence (Proposition \ref{ks-modif}), cf. \cite[Proposition 1.7]{Pal90}: $\xi$ being semi-universal, $\theta_\xi$ is injective, thus we may again identify $T^0(X_0,\mf{m})$ and $T^0(\xi,(\C,\mf{m}))$. Together with the surjectivity of the map $ev'$ (proved in Lemma \ref{evprime-surj}) we obtain that $ev'$ induces a surjective $\OO_S$-linear map 
\[T^0(\xi,(\OO_S,J_\sigma))/T^0(X/S,J_\sigma)\to T^0_\bullet(X_0,\mf{m}).\] 
The $\OO_S$-module on the left is isomorphic to the submodule $Im(T^0(\xi,(\OO_S,J_\sigma))\to T^0(S))=Ker(\Theta_\xi)$ of the finitely generated $\OO_S$-module $T^0(S)$, hence it is itself finitely generated. Now applying $-\otimes_{\OO_S} \C$ yields a surjective map
\[T^0(\xi,(\OO_S,J_\sigma))/T^0(X/S,J_\sigma)\otimes_{\OO_S} \C \to T^0_\bullet(X_0,\mf{m})\otimes_{\OO_S} \C\simeq T^0_\bullet(X_0,\mf{m}),\]
whose source is a finite-dimensional $\C$-vector space, hence $\dim_\C T^0_\bullet(X_0,\mf{m})<\infty$, too.

We now prove $(2)$, which is done analogously to \cite[Proposition 3.1]{Pal90}. $\eps$ is well-defined: Let $X$ be defined by $F\in \OO_S\{\x\}$, in addition we may assume that $\sigma$ is the zero section, i. e. $J_\sigma=(\x)\OO_S\{\x\}$. If $\delta\in T^0(X_0,\mf{m})$ lifts to $\tilde\delta\in T^0(X/S,J_\sigma)$, then there exist $G_i\in (\x)\OO_S\{x\}$ and $H\in \OO_S\{\x\}$ such that $\tilde\delta$ is the class of $\sum_i G_i\frac{\p}{\p x_i}$ and $\sum_i G_i\frac{\p F}{\p x_i}=H\cdot F$. Since $F$ and its partial derivatives with respect to the $x_i$ form a regular sequence in $\OO_S\{\x\}$, this implies that $H$ is in the ideal generated by $\frac{\p F}{\p x_1},\dots,\frac{\p F}{\p x_n}$, i. e. $h=H_{|s=0}\in J(f)$. 

It remains to verify the exactness of the sequence, which is obvious at the two terms on the right. If $h \bmod{J(f)}\in Ker(m_f)$, then $hf=\sum_i g_i\frac{\p f}{\p x_i}$ with $g_i\in (\x)$, i. e. $h=\eps(\delta)$ for the derivation $\delta:\OO_{X_0}\to \mf{m}$ induced by $\sum_i g_i\frac{\p}{\p x_i}$. Since the inclusion $Im(\eps)\subseteq Ker(m_f)$ is clear, it is only left to show that $\eps$ has trivial kernel. But if $\delta$ is induced by a derivation $\delta'$ of $\OO_n$ with $\sum_i g_i\frac{\p f}{\p x_i}=h\cdot f$ and $h=\sum_i h_i\frac{\p f}{\p x_i}$, then its class equals the class of $\sum_i (g_i-f\cdot h_i)\frac{\p}{\p x_i}$, so we may assume that for its lift $\delta'$ we have $\delta'(f)=0$. Thus its coefficients are a syzygy of the partial derivatives of $f$ which form a regular sequence in $\OO_n$. Hence $\delta$ belongs to the submodule generated by the classes of the derivations of the form $\frac{\p f}{\p x_j}\frac{\p}{\p x_i}-\frac{\p f}{\p x_i}\frac{\p}{\p x_j}$, which clearly can be lifted to $T^0(X/S)$.

The last claim is now a consequence of Lemma \ref{dimT1sec}, which implies that
\[\dim_C \OO_n/(f,(x)J(f))=\dim_C T^1(X_0,\mf{m})+1 = \tau(X_0)+n.\]
Now let, in addition, $\mu(X_0):=\dim_C \OO_n/J(f)$ be the Milnor number of $X_0$. From the fact that the partial derivates of $f$ form a regular sequence, we conclude that they are linearly independent modulo $(x)J(f)$, hence $\dim_C \OO_n/(x)J(f)=\mu+n$. Altogether we calculate using (\ref{T0red-seq}):
\begin{align*}
\dim_C T^0_\bullet(X_0) & = \dim_C \OO_n/J(f)-\dim_C \OO_n/(x)J(f) + \dim_C \OO_n(f,(x)J(f))\\
                        & = \mu - (\mu + n) + (\tau(X_0)+n) \\
			& = \tau(X_0),
\end{align*}			
which proves the last assertion.
\end{proof}

\subsection{Quasihomogeneous complete intersections}

In general, the modular stratum of an isolated singularity $X_0$ carries a non-reduced structure. This is, for instance, the case for semi-quasihomogeneous hypersurface singularities. However, it is proved in \cite[section 6.2]{Alex85}, that the modular stratum is reduced and smooth in case $X_0$ is a quasihomogeneous isolated complete intersection singularity. As an application of the criteria for modularity established above, we extend this result to defomations with section:

\begin{proposition}\label{qh-modstratum}
The modular stratum with respect to the functor $Def^s_{X_0}$ of a quasihomogeneous isolated complete intersection singularity $X_0$ is reduced and smooth. 
\end{proposition}

\begin{proof}
Let $X_0$ be defined by quasihomogeneous polynomials $f=(f_1,\dots,f_k)\in\OO_n^k$ of degrees $d_1,\dots,d_k$ with respect to some positive integer weights $w_1,\dots,w_n$. For any $1\le i_1<\ldots<i_{k+1}\le n$ one obtains a {\em Hamiltonian derivation} by cofactor expansion with respect to the first line of the symbolic matrix
\[H(i_1,\dots,i_{k+1}):=\det\left|\BMA
\frac{\p}{\p x_{i_1}} &\cdots & \frac{\p}{\p x_{i_{k+1}}} \\
\frac{\p f_1}{\p x_{i_1}} &\cdots & \frac{\p f_1}{\p x_{i_{k+1}}} \\
\vdots & & \vdots \\
\frac{\p f_k}{\p x_{i_1}} &\cdots & \frac{\p f_k}{\p x_{i_{k+1}}} 
\BME\right|,\]
and it is proved in \cite[section 6.1]{Alex85} that the Hamiltonian derivations together with the Euler derivation $\delta_E:=\sum_{i=1}^n w_i x_i\frac{\p}{\p x_i}$ generate the $\OO_{X_0}$-module $T^0(X_0)=\Der_\C(\OO_{X_0})$. Note that in particular all these derivations take their image in the maximal ideal $\mf{m}$ of $\OO_{X_0}$, so in this situation we have $T^0(X_0)=T^0(X_0,\mf{m})$ and we can take them as generators of the latter module, too.

Now let $\xi:X\to S$ be a semi-universal deformation with section $\sigma$ of $X_0$. In particular, by the results of section 2.3, $S\simeq(\C^{\tau^s},0)$ is smooth. Obviously, all $H(i_1,\dots,i_{k+1})$ can be lifted to $T^0(X/S,J_\sigma)$, so we only have to consider $\delta_E$. The bracket $[-,-]$ in cotangent cohomology induces an element $[\delta_E,-]\in \mathrm{End}_\C(T^1(X_0,\mf{m}))$ that gives a decomposition $T^1(X_0,\mf{m})=\bigoplus_{\nu\in \Z}T^1(X_0,\mf{m})_\nu$ into a direct sum of its eigenspaces $T^1(X_0,\mf{m})_\nu$ with respect to the eigenvalue $\nu$. Using the description of the tangent space to the modular stratum $M$ in Corollary \ref{tangmod} we obtain $T(M)\simeq T^1(X_0,\mf{m})_0$. 

But now it is clear that there are no further obstructions to lifting $\delta_E$: If $F=f+\sum_{i=1}^{\tau^s}s_i g^{(i)}$ is a semi-universal family with monomials $g^{(i)}$ as in section 2.3 above, then we may assume that $g^{(1)},\dots,g^{(r)}$ represent a basis of $T^1(X_0,\mf{m})_0$. From what was said above, it then follows that $T(M)\simeq ((s_1,\dots,s_r)/(s_1,\dots,s_r)^2)^*$, and we can lift $\delta_E$ to the restriction of $\xi$ to the smooth subspace $(\C^r\times\{0\},0)$ of $S$, i. e. $\xi_{|(\C^r\times\{0\},0)}$ is the maximal modular deformation inside $\xi$.
\end{proof}

The argument given above is entirely the same in the context of deformation without section, it is a consequence of the identity $T^0(X_0)= T^0(X_0,\mf{m})$ for singularities of this type. As a corollary we obtain:

\begin{corollary}
For any quasihomogeneous isolated complete intersection singularity $X_0$ the modular strata with respect to both deformation functors $Def_{X_0}$ and $Def^s_{X_0}$ coincide.
\end{corollary}

%--- Flatness and modular deformations ----------------------------------------

\section{Flatness and modular deformations}

\subsection{Modularity as flatness of the first Tjurina module}

Now choose a minimal embedding $(X_0,0)\subseteq(\C^n,0)$. Let $X_0$ be defined by $f_1,\dots,f_k\in\OO_n$ and let $\xi:X\to S$ be a deformation of $X_0$, given by $F_1,\dots,F_k\in \OO_S\{\x\}$. Denote by $J(f)$ the Jacobian matrix of $f_1,\dots,f_k$ and by $J(F)$ the relative (with respect to $x_1,\dots,x_n$) Jacobian matrix of $F_1,\dots,F_k$. Finally, set $\tilde T^1(X/S):=\OO^k_X/J(F)$. With these notions we can interpret the criterion in Theorem \ref{modcrit} as flatness of the module $\tilde T^1(X/S)$:

\begin{proposition}\label{flatcrit}
Let $M\subseteq S$ be a subspace in the base of a semi-universal deformation $\xi:X\to S$ of unobstructed singularity $X_0$. Then $M$ is modular if and only of $\tilde T^1(X/S)\otimes_{\OO_S}\OO_M$ is a flat $\OO_M$-module. 
\end{proposition}

\begin{proof}
There is the following commutative diagram with exact rows:
\[\xymatrix{
0\ar[r] & T^0(X/S)\ar[r]\ar[d]^{ev} & \OO^n_X\ar[d]\ar[r]^{J(F)} & \OO^k_X\ar[d]\ar[r] & \tilde T^1(X/S)\ar[d]\ar[r] & 0 \\
0\ar[r] & T^0(X_0)\ar[r] & \OO^n_{X_0}\ar[r]^{J(f)} & \OO^k_{X_0}\ar[r] & T^1(X_0)\ar[r] & 0
}\]
Since $X_0$ was assumed to be unobstructed, $\OO_X$ is smooth and it follows that it is a free $\OO_S$-module. Hence we can interpret $T^0(X/S)$ (resp. $T^0(X_0)$) as the syzygy module of the columns of the presentation matrix of $\tilde T^1(X/S)$ (resp. $T^1(X_0)$). Thus, by the lifting criterion for flatness (see \cite[§ 7]{Mat86}, for example),  $\tilde T^1(X/S)\otimes_{\OO_S}\OO_M=\tilde T^1(X_{|M}/M)$ is a flat $\OO_M$-module if and only if $ev_{|M}$ is surjective. 
\end{proof}

If $X_0$ is a complete intersection then $\tilde T^1(X/S)$ equals $T^1(X/S)$, hence:

\begin{corollary}\label{modstrat}
The modular stratum of an isolated complete intersection singularity $X_0$ equals the flattening stratum of the relative Tjurina module $T^1(X/S)$, where $\xi:X\to S$ is a semi-universal deformation of $X_0$.
\end{corollary}

\begin{remark}
This result is extended to reduced space curve singularities in \cite{Mar03} using their determinantal structure (defining equations of any deformation of such a singularity are obtained as maximal minors of some $q\times(q-1)$-matrix $A$, and the relativ normal module of the deformation has a presentation matrix  whose entries are the $(q-2)$-minors of $A$). However it is not clear  whether these assertions can be carried over to arbitrary isolated singularities, or at least to arbitrary unobstructed ones.
\end{remark}

Nevertheless, analogous statements hold for deformations with section, which we are now going to derive:

\begin{proposition}\label{flatcrit-sec}
Let $\xi:X\to S$ be a semi-universal deformation of an ubstructed singularity $X_0$ with section $\sigma:S\to X$. Let $X_0$ be defined by $f_1,\dots,f_k\in \C\{x_1,\dots,x_n\}$, $X$ be given by $F_1,\dots,F_k\in\OO_S\{x_1,\dots,x_n\}$. Denote $J_\sigma:=Ker(\sigma^*)\subseteq \OO_X$, and let $J(F):=\left(\frac{\p F_i}{\p x_j}\right)_{i,j}\in Mat(k,n;\OO_X)$ be the relative Jacobian matrix of $F_1,\dots,F_k$.

Then a subgerm $M\subseteq S$ is modular if and only if $\tilde T^1(X/S,J_\sigma)\otimes_{\OO_S}\OO_M$ is a flat $\OO_M$-module, where
\[\tilde T^1(X/S,J_\sigma):=\OO^k_{X}/J_\sigma\cdot J(F).\]
\end{proposition}

\begin{proof}
We may assume that $J_\sigma=(\x)\OO_X$. There is the following commutative diagram whose rows are exact sequences:
\[\xymatrix{
0\ar[r]& T^0(X/S,J_\sigma)\ar[r]\ar[d]_{ev} & \OO_X^{n^2} \ar[rrr]^{(x_1 J(F),\dots,x_nJ (F))}\ar[d] &&& \OO_X^k\ar[d]\ar[r]& \tilde T^1(X/S,J_\sigma)\ar[d]\ar[r]& 0 \\
0\ar[r]& T^0(X_0,\mf{m})\ar[r] & \OO_{X_0}^{n^2}\ar[rrr]^{(x_1 J(f),\dots,x_n J(f))} &&& \OO_{X_0}^k\ar[r] & \OO_{X_0}/(\mf{m}J(f))\ar[r] & 0.
}\]
As in the proof of Proposition \ref{flatcrit} we can interpret $T^0(X/S,J_\sigma)$ (resp. $T^0(X,\mf{m})$) as syzygy modules of the modules $\tilde T^1(X/S,J_\sigma)$ (resp. $\OO_{X_0}/(\mf{m}J(f))$). So the surjectivity of $ev_{|M}$ means precisely that every syzygy of $\tilde T^1(X/S,J_\sigma)$ over the special fibre $s=0$ lifts to a syzygy over $M$, i. e. $\tilde T^1(X/S,J_\sigma)\otimes_{\OO_S}\OO_M=\tilde T^1(X_{|M}/M,J_{\sigma|M})$ is $\OO_M$-flat.
\end{proof}

\begin{corollary}\label{modstrat-sec}
If $X_0$ is an isolated complete intersection, then the modular stratum (with section) of $X_0$ coincides with the flattening stratum of the relative Tjurina module $T^1(X/S,J_\sigma)$, where $S$ is the base space of a semi-universal deformation $X\to S$ of $X_0$ and $J_\sigma=ker(\sigma^*)\subseteq\OO_X$ corresponds to the section $\sigma:S\to X$.
\end{corollary}

\begin{proof}
$X_0$ being a complete intersection we have, for any subspace $M\subseteq S$, an exact sequence
\[0\to T^1(X/S,J_\sigma)_{|M}\to \tilde T^1(X/S,J_\sigma)_{|M}\stackrel{r}{\to} \OO_M\to 0\]
of $\OO_M$-modules, where $r$ is the residue map. Hence $T^1(X/S,J_\sigma)_{|M}$ is $\OO_M$-flat if and only if $\tilde T^1(X/S,J_\sigma)_{|M}$ is.
\end{proof}

These characterizations of modularity as flatness of a suitable $\OO_S$-module make it possible to compute the modular stratum of $X_0$, at least up to a given order. Computing the constant rank stratum or the annihilator of the associated torsion module whose vanishing implies flatness (see e.g. \cite[section 7.3]{GP02}) is not a feasible approach in this situation: The objects involved are too complicated for computations. Instead, one can determine, starting with the base field $\C$, the maximal small extension in $\OO_S$ that preserves flatness. The resulting algorithm and details on its implementation in the \Singular-library \verb+modular.lib+ (\cite{modular}) are described in \cite{Mar02}, where also several non-trivial examples of modular strata for the functor $Def_{X_0}$ are presented. 

\begin{remark}
Of course one may consider the flattening stratum of $T^1(X/S)$ for any deformation $\xi:X\to S$ of $X_0$. But if $\xi$ is not a monodeformation, there is no connection to modular subspaces: The latter must be trivial, whereas the former may not, as the following example shows. Take an arbitrary hypersurface singularity $X_0$ defined by $f\in \OO_n$ and $g\in (f,J(f))$. Then $f+\eps g$ gives a trivial deformation $X_D$ of $X_0$ over the double point $D$, its Kodaira-Spencer map is the zero map. On the other hand $T^1(X_D/D)$ is a flat $\OO_D=\C[\eps]$-module: 

Indeed, let $rf+\sum_k r_k \frac{\p f}{\p x_k}=0$ be a relation of $f,J(f)$. Then $r(f+\eps g)+\sum_k r_k\frac{\p}{\p x_k}(f+\eps g)=\eps(rg+\sum_k r_k\frac{\p g}{\p x_k})$, and one checks that this again belongs to $\eps(f,J(f))$, so is equal to $\eps(pf+\sum_k p_k\frac{\p f}{\p x_k})$ for some $p,p_k\in \OO_n$. Then $(r-\eps p,r_1-\eps p_1,\dots,r_n-\eps p_n)$ lifts the given relation to $\OO_n[\eps]$. Hence any trivial first-order infinitesimal deformation of a hypersurface has flat relative $T^1(X_D/D)$.
\end{remark}

\subsection{Comparing flatness conditions}
As established above, flatness of the module $T^1(X/S)$ (resp. $T^1(X/S,J_\sigma)$ has an interpretation as modularity of the corresponding deformation $\xi:X\to S$ (with section $\sigma$). One can then ask about the relationship to flatness of the $0$-th cotangent cohomology module. The result is:

\begin{proposition}
Let $\xi:X\to S$ be a deformation of an isolated complete intersection singularity $X_0$. If $T^1(X/S)_{|M}$ is a flat $\OO_M$-module, then so is $T^0(X/S)_{|M}$, i. e. the flattening stratum of $T^1(X/S)$ is contained in the flattening stratum of $T^0(X/S)$, and this inclusion may be strict.\\
The analogous statement holds for deformations with section.
\end{proposition}

\begin{proof}
By splitting the short exact sequence
\[0\to T^0(X/S)_{|M}\to \OO_{X_{|M}}^n \to \OO_{X_{|M}}^k \to T^1(X/S)_{|M}\simeq \tilde T^1(X_{|M}/M)\to 0\]
from the proof of Proposition \ref{flatcrit} into two short exact sequences we deduce that $T^0(X/S)_{|M}$ is flat if $T^1(X/S)_{|M}$ is. In the case of a deformation with section we use the sequence in the proof of Proposition \ref{flatcrit-sec} instead.

On the other hand, let $X_0$ be a quasihomogeneous complete intersection defined by $f=(f_1,\dots,f_k)\in\OO_n$ and take a basis monomial of $T^1(X_0)$ resp. $T^1(X_0,\mf{m})$ which is of different weighted degree than $f$. Let $\xi:X_D\to D$ be the deformation over the double point $D$ defined by $f+\eps g$ (with the obvious section $\sigma$). From Proposition \ref{qh-modstratum} we conclude that this monodeformation has trivial modular stratum, hence the flattening stratum of $T^1(X_D/D)$ resp. $T^1(X_D/D,J_\sigma)$ is trivial. But, as a submodule of $\OO_{X_D}^n$ (resp. $\OO_{X_D}^{n^2}$), $T^0(X_D/D)$ (resp. $T^0(X_D/D,J_\sigma)$) is a torsion-free, hence a flat $\OO_D=\C[\eps]$-module.
\end{proof}

We finish by comparing the modular strata of a hypersurface with respect to both deformation functors under consideration. First of all note the following fact:

\begin{lemma}\label{mod-mods}
Let $\xi:X\to S$ be a deformation of $X_0$ with section $\sigma:S\to X$. If $M\subseteq S$ is modular with respect to the functor $Def_{X_0}$ (i. e. considering $\xi$ as an 'ordinary' deformation), then it is modular with respect to $Def^s_{X_0}$, too.
\end{lemma}

\begin{proof}
Suppose $M$ is $Def_{X_0}$-modular. Let $\ph:T\to M$ and $\psi:T\to S$ be morphisms such that $\ph^*(\xi_{|M})$ and $\psi^*(\xi)$ are isomorphic deformations with section of $X_0$. Then they also give the same element of $Def_{X_0}(S)$, hence $\ph=\psi$ by assumption. Thus $M$ is $Def^s_{X_0}$-modular, too.
\end{proof}

We can prove a converse of this statement if we restrict ourselves to deformations with {\em singular section} of $X_0$: If $X_0$ is equidimensional of dimension $d$, then a deformation with singular section of $X_0$ is a deformation $\xi:X\to S$ with section $\sigma:S\to X$, such that there is a factorization
\[\xymatrix{
S\ar[rr]^\sigma\ar@{.>}[dr] &&X, \\
& V(F^d(\Omega_{X/S}))\ar@{^{(}->}[ur]}\]
where $F^d(\Omega_{X/S})$ is the $d$-th Fitting ideal of the module of Kähler differentials of the $\OO_S$-module $\OO_X$, $C_\xi:=V(F^d(\Omega_{X/S}))$ is the {\em critical locus} of $\xi$. One then has the following properties:
\begin{enumerate}
\item This definition gives rise to a subfunctor $Def^{ss}_{X_0}$ of $Def^s_{X_0}$. A versal deformation of $X_0$ with respect to $Def^{ss}_{X_0}$ is obtained as follows: Let $X\to S$ be a versal deformation of $X_0$ (with respect to $Def_{X_0})$, then, as seen in section 2, $X\times_S X\to X$ is versal for $Def^s_{X_0}$, and then the deformation induced by the base change $V(F^d(\Omega_{X/S}))\hookrightarrow X$ is a versal deformation with singular section of $X_0$ (cf. \cite[2.5.1]{Buch81}).
\item If $f\in \OO_n$ defines an isolated hypersurface singularity, then one can easily define a semi-universal deformation with singular section of $X_0:=(V(f),0)$ by $F:=f+\sum_i s_i g_i$, where $\{g_i\}_i$ represent a $\C$-basis of $\mf{m}\cdot T^1(X_0,\mf{m})$, $\mf{m}$ the maximal ideal of $\OO_{X_0}$.
\end{enumerate}

The latter construction cannot be generalized to singularities of higher codimension. For hypersurfaces, we can interpret deformations with singular section as deformations leaving the embedding dimension constant, the terms of its defining equation then belong $\mf{m}^2$. This is why we obtain the following result only for hypersurfaces.

\begin{proposition}\label{flat-ss}
Let $\xi:X\to S$ be a deformation with singular section of a hypersurface $X_0$ and let $M\subseteq S$. If $T^1(X/S)_{|M}$ is a flat $\OO_M$-module, then so is $T^1(X/S,J_\sigma)_{|M}$.
\end{proposition}

\begin{proof}
The short exact sequence of $\OO_{X_{|M}}$-modules
\[0\to J_{\sigma|M} \to \OO_{X_{|M}} \to \OO_{X_{|M}}/J_{\sigma|M}\simeq\OO_M\to 0\]
yields the long exact cohomology sequence
\begin{multline}\label{flat-ss-proof}
\xymatrix{
0\ar[r]& T^0(X/S,J_\sigma)_{|M}\ar[r]^(0.55){\iota_0}&  T^0(X/S)_{|M}\ar[r]^(0.45){\pi_0}&  T^0(X/S,\OO_S)_{|M}\ar[r]^(0.75){\delta_0}& } \\
\xymatrix{
\ar[r]^(0.25){\delta_0} & T^1(X/S,J_\sigma)_{|M} \ar[r]^(0.55){\iota_1} & T^1(X/S)_{|M}\ar[r]^(0.45){\pi_1} & T^1(X/S,\OO_S)_{|M}\ar[r] & 0.}
\end{multline}
Let $f\in\OO_n$ define $X_0$. Since $X_{|M}\to M$ is a deformation with singular section of $X_0$, we may assume that $X_{|M}$ is defined by $F\in (\x)^2\OO_M\{\x\}$. In this situation we have $T^1(X/S,\OO_S)_{|M}\simeq\OO_M$ since it is the cokernel of
\begin{align*}
\Der_{\OO_M}(\OO_M\{\x\},\OO_M) & \to \Hom_{\OO_{X|M}}((F)/(F^2),\OO_M)\simeq \OO_M,\\
\delta & \mapsto \delta(F)
\end{align*}
which is the zero map as $F\in(\x)^2\OO_M\{\x\}$. In addition, $Ker(\iota_1)$ is the $\OO_M$-module generated by $\frac{\p F}{\p x_1},\dots,\frac{\p F}{\p x_n}$ modulo $(F)+(x_j\frac{\p F}{\p x_i})_{i,j}$. But this is a free $\OO_M$-module: namely, suppose $\sum_i G_i(s)\frac{\p F}{\p x_i}= H\cdot F+\sum_{i,j} H_{ij}x_j \frac{\p F}{\p x_i}$. Since $(F,\frac{\p F}{\p x_1},\dots,\frac{\p F}{\p x_n})$ form a regular sequence, this implies that, for all $i$, $G_i(s)-\sum_j H_{ij}x_j\in (F,\frac{\p F}{\p x_1},\dots,\frac{\p F}{\p x_n})$, i. e.
\[G_i(s)= \sum\nolimits_j H_{ij}x_j + A\cdot F + \sum\nolimits_k B_k \tfrac{\p F}{\p x_k}\]
for some $A,B_k\in \OO_M\{\x\}$. But this is only possible if $G_i=0$ since $F,\frac{\p F}{\p x_1},\dots,\frac{\p F}{\p x_n}\in (\x)\OO_M\{\x\}$.

Thus the long exact sequence (\ref{flat-ss-proof}) can be split into two short exact sequences
\[0\to Ker(\iota_1)\simeq \OO_M \to T^1(X/S,J_\sigma)_{|M} \to Im(\iota_1)\to 0\]
and
\[0\to Ker(\pi_1)\to T^1(X/S)_{|M} \to \OO_M\to 0,\]
from which we derive that if $T^1(X/S)_{|M}$ is $\OO_M$-flat, then so is $T^1(X/S,J_\sigma)_{|M}$.
\end{proof}

Combining Lemma \ref{mod-mods} and Proposition \ref{flat-ss} with the modularity criteria Corollary \ref{modstrat} and Corollary \ref{modstrat-sec}, we obtain:

\begin{corollary}
Let $\xi:X\to S$ be a deformation with singular section $\sigma:S\to X$ of a hypersurface $X_0$. Then 
\begin{enumerate}
\item The flattening strata of $T^1(X/S)$ and $T^1(X/S,J_\sigma)$ coincide.
\item If $\xi$ is a monodeformation for both deformation functors $Def_{X_0}$ and $Def^s_{X_0}$, then the maximal modular subspaces inside $S$ with respect to these two functors coincide.
\end{enumerate}
\end{corollary}

%-- Literatur -----------------------------------------------------------------

\newpage

\addcontentsline{toc}{section}{References}
\bibliographystyle{amsalpha}
\bibliography{../literatur}

\end{document}